\documentclass[12pt]{article}

\usepackage{mathtools}
\usepackage{amssymb}
\usepackage{amsthm}
\usepackage{tikz}
\usepackage{fullpage}
\usepackage{url}

\DeclareSymbolFont{bbold}{U}{bbold}{m}{n}
\DeclareSymbolFontAlphabet{\mathbbold}{bbold}

\newcommand{\N}{\mathbb{N}}
\newcommand{\Z}{\mathbb{Z}}

\newcommand{\R}{\mathbb{R}}
\newcommand{\C}{\mathbb{C}}

\newcommand{\mm}{\mathrm{m}}
\newcommand{\ii}{\mathrm{i}}
\newcommand{\dd}{\mathrm{d}}
\newcommand{\ee}{\mathrm{e}}
\newcommand{\nlH}{\mathrm{nlH}}

\newcommand{\w}{\mathrm{w}}
\renewcommand{\hom}{\mathrm{hom}}

\renewcommand{\epsilon}{\varepsilon}

\DeclareMathAccent{\Circ}{\mathalpha}{operators}{"17}

\newcommand{\Gto}{\stackrel{G}{\to}}
\newcommand{\Hto}{\stackrel{H}{\to}}



\DeclareMathOperator{\sym}{sym}
\DeclareMathOperator{\dive}{div}
\DeclareMathOperator{\curl}{curl}
\DeclareMathOperator{\grad}{grad}

\DeclareMathOperator{\kar}{ker}
\DeclareMathOperator{\rge}{ran}
\DeclareMathOperator{\dom}{dom}

\renewcommand{\Re}{\operatorname{Re}}
\renewcommand{\hat}{\widehat}

\renewcommand{\d}{\,\mathrm{d}}

\let\phi\varphi

\let\geq\geqslant

\makeatletter
\def\@row#1,{#1\@ifnextchar;{\@gobble}{&\@row}}
\def\@matrix{%
    \expandafter\@row\my@arg,;%
    \@ifnextchar({\\ \get@in@paren{\@matrix}}{\after@matrix}%
    }
\def\matrixtype#1#2#3{%
    \ifmmode\def\after@matrix{\end{#2}\right#3}%
    \else\def\after@matrix{\end{#2}\right#3$}$\fi\iffalse$\fi
    \left#1\begin{#2}\get@in@paren{\@matrix}%
    }
\def\@column#1,{#1\@ifnextchar;{\@gobble}{\\ \@column}}
\newcommand\vect{}
\def\svect(#1){\left(\begin{smallmatrix}\@column#1,;\end{smallmatrix}\right)}
\def\vect{\get@in@paren{\@vect}}
\def\@vect{\left(\begin{matrix}\expandafter\@column\my@arg,;\end{matrix}\right)}
\def\get@in@paren#1({\def\my@arg{}\def\my@rest{}\def\after@get{#1}\get@arg}
\let\e@a\expandafter
\def\get@arg#1){\e@a\kl@test\my@rest#1(;}
\def\kl@test#1(#2;{\e@a\def\e@a\my@arg\e@a{\my@arg#1}%
                   \ifx:#2:\let\my@exec\after@get
                   \else\let\my@exec\get@arg
                        \e@a\def\e@a\my@arg\e@a{\my@arg(}%
                        \def@rest#2;%
                   \fi\my@exec}
\def\def@rest#1(;{\def\my@rest{#1\kl@zu}}
\def\kl@zu{)}

\makeatletter
\newcommand\MyPairedDelimiter{%
  \@ifstar{\My@Paired@Delimiter{{}}}
          {\My@Paired@Delimiter{}}%
}
\newcommand\My@Paired@Delimiter[4]{%
  \newcommand#2{%
    \@ifstar{\start@PD{#1}{\delimitershortfall=-1sp}{#3}{#4}}
            {\start@PD{#1}{}{#3}{#4}}%
  }%
}
\newcommand\start@PD[5]{%
  #1\mathopen{\mathpalette\put@delim@helper{\put@delim{#2}{#3}{.}{#5}}}%
  #5%
  \mathclose{\mathpalette\put@delim@helper{\put@delim{#2}{.}{#4}{#5}}}%
}
\newcommand\put@delim@helper[2]{%
  \hbox{$\m@th\nulldelimiterspace=0pt #2#1$}%
}
\newcommand\put@delim[5]{%
  \setbox\z@\hbox{$\m@th#5{#4}$}%
  \setbox\tw@\null
  \ht\tw@\ht\z@ \dp\tw@\dp\z@
  #1#5%
  \left#2\box\tw@\right#3%
}

\makeatother
\MyPairedDelimiter*{\abs}{\lvert}{\rvert}
\MyPairedDelimiter*{\norm}{\lVert}{\rVert}
\MyPairedDelimiter{\set}{\{}{\}}

\theoremstyle{plain} % default
\newtheorem{theorem}{Theorem}[section]
\newtheorem{corollary}[theorem]{Corollary}

\newtheorem{proposition}[theorem]{Proposition}
\theoremstyle{definition}
\newtheorem{example}[theorem]{Example}

\newtheorem*{definition}{Definition}
\newtheorem{remark}[theorem]{Remark}

\usepackage{enumitem}

\setenumerate[1]{nolistsep} % Only the level 1
\setenumerate[2]{nolistsep} % Only the level 2

\begin{document}

\title{Homogenisation and the Weak Operator Topology}

\author{Marcus Waurick}

\date{}

\maketitle

\begin{abstract} This article surveys results that relate homogenisation problems for partial differential equations and convergence in the weak operator topology of a suitable choice of linear operators.  More precisely, well-known notions like $G$-convergence, $H$-convergence as well as the recent notion of nonlocal $H$-convergence are discussed and characterised by certain convergence statements under the weak operator topology. Having introduced and described these notions predominantly made for static or variational type problems, we further study these convergences in the context of dynamic equations like the heat equation, the wave equation or Maxwell's equations.  The survey is intended to clarify the ideas and highlight the operator theoretic aspects of homogenisation theory in the autonomous case. 
\end{abstract}

Keywords: homogenisation, evolutionary equations, $G$-convergence, $H$-convergence, nonlocal $H$-convergence, heat conduction, wave equation, Maxwell's equations

MSC 2010: Primary 35B27 74Q10 74Q05; Secondary  35L04  35Q61

\medmuskip=4mu plus 2mu minus 3mu
\thickmuskip=5mu plus 3mu minus 1mu
\belowdisplayshortskip=9pt plus 3pt minus 5pt

\section{Introduction}\label{s:int}

The theory of homogenisation is concerned with the analysis of equations with highly oscillatory coefficients and whether these coefficients can be `averaged out' in a suitable sense. In the classical setting one considers an equation with periodic coefficients and asks oneself what happens if the period length tends to $0$. We refer to some standard literature introducing the topic: \cite{Bensoussan1978,Cioranescu1999,Tartar2009}. 

In this survey, we draw the relation of this study to convergence in the weak operator topology. Note that in a very general context, this connection is well-known, see e.g.~\cite[Lemma 6.2]{Tartar2009} or \cite{Zhikov1979,Spagnolo1967,Spagnolo1976}. In the present article, we want to provide a more detailed perspective. In fact, we will characterise classical notions in (non-periodic) homogenisation theory by suitable operators converging in the weak operator topology. We shall highlight the most important results obtained by the author starting in \cite{W11_P} and culminating in the rather recent contribution \cite{W18_NHC} in the context of homogenisation theory of partial differential equations. For ordinary differential equations, we shall refer to \cite{W12_HO,W14_G} for the original research papers as well as to \cite{W16_H} for a round up of these results in the up-to-date most general context of so-called evolutionary mappings. We shall also refer to the references highlighted there for an account of results in the literature. We just briefly mention that it is possible to detour the concept of so-called Young measures in the context of homogenisation of ordinary differential equations as it has been used in \cite[Chapter 23]{Tartar2009} and to obtain equally precise results.

Before we start to provide some more details, let us mention the most important conceptual difference to other approaches in the literature. In classical homogenisation theory particularly in time-dependent problems, the coefficients of the spatial derivative operators are assumed to be highly oscillatory. Rewriting the equation in question, the philosophy is rather to consider spatial derivative operators that have no coefficients multiplied to them. This allows to have a constant domain of definition for the operators considered. Moreover, with this perspective it is possible to straightforwardly address equations of mixed type (also highly oscillatory versions of it) in a comprehensive space-time setting. Also with this methodology it is easier to spot the occurrence of memory effects due to the homogenisation process, see e.g.~\cite{W14_FE}.

Focussing on partial differential equations in the following, we begin with a motivation to relate  highly oscillatory functions and the weak operator topology. Take a $1$-periodic function $a\in L^\infty(\R)$ and consider $a_n\colon x\mapsto a(nx)$ for all $n\in \N$. Note that the period length of $a_n$ is bounded above by $1/n$. It is then not very difficult to see that for all bounded intervals $I\subseteq \R$, we have, as $n\to\infty$,
\[
    \int_I a_n(x) \d x\to \int_I \mathfrak{M}(a) \d x,
\]
where $\mathfrak{M}(a)\coloneqq \int_0^1 a(x) \d x$, see \cite[Theorem 2.6]{Cioranescu1999}. Using linearity of the integral and that simple functions are dense in $L^1(\R)$, we deduce that $a_n \to \mathfrak{M}(a)$ in $\sigma(L^\infty,L^1)$ as $n\to\infty$. By Cauchy--Schwarz' inequality, we also see that any $L^1$-function can be written as the product of two $L^2$-functions. Thus, we deduce that
\[
    a_n \to \mathfrak{M}(a) \in \mathcal{B}_\textnormal{w}(L^2(\R))\quad(n\to\infty),
\]
where $\mathcal{B}_\textnormal{w}(L^2(\R))$ is the set of bounded linear operators on $L^2(\R)$ endowed with the weak operator topology.

We shall describe the contents of this survey in more detail next. We will head off with local and afterwards nonlocal coefficients for elliptic partial differential equations. After a small interlude on how to obtain a Hilbert space setting for time-dependent problems, we shall discuss  homogenisation theorems of partial differential equations with a spatial operator that has compact resolvent. We conclude this survey with a result combining all the results obtained before. The results will always be formulated in the way that certain solution operators fo the equations in question converge in a suitable weak operator topology. We will also mention, where the convergence might even be improved. 

Although we will treat homogenisation that go well beyond the periodic case, the reader might think of the periodic case as a first particular application of the abstract results. This survey focusses on the results, we will largely refrain from providing the respective proofs, but we shall sketch some applications. In order to present a reasonably widespread amount of contents we keep the applications and examples on a mere informal level and refer to the original papers for the details.

Throughout this manuscript Hilbert spaces are anti-linear in the first and linear in the second component. The arrow $\rightharpoonup$ symbolises weak convergence and $\mathcal{B}(H)$ denotes the set of bounded linear operators on a Hilbert space $H$, we write $\mathcal{B}_\w(H)$ if we want to stress that $\mathcal{B}(H)$ is endowed with the weak operator topology. As usual we write $\mathcal{B}(H,K)$ for operators acting from $H$ to Hilbert space $K$. The letters $\alpha$ and $\beta$ will always denote strictly positive real numbers with $\alpha<\beta$. The $*$ as an upper index denotes the operator adjoint as well as complex conjugate.

\section{Time-independent problems -- local coefficients}

\emph{Periodic problems.} Historically, one of the first examples in the study of homogenisation theory was concerned with periodic elliptic problems. More precisely, let $a\in L^\infty(\R^N)^{N\times N}$ satisfy for almost all $x\in \R^N$ and all $k\in \Z^N$
\[
  \frac12\left(  a(x)+a(x)^*\right)\geq \alpha\text{ and } a(x)=a(x+k),
\]
where the first inequality holds in the sense of positive definiteness for some $\alpha>0$. We set $a_n\coloneqq a(n\cdot)$. With these settings, we can state the by now classical theorem on elliptic homogenisation problems;

\begin{theorem}[{{classical, see e.g.~\cite[Theorem 6.1]{Cioranescu1999}}}] Let $\Omega\subseteq \R^N$ open, bounded. Let $f\in H^{-1}(\Omega)$ and $u_n\in H_0^1(\Omega)$ be such that 
\[
            \langle a_n \grad u_n , \grad \phi\rangle = f(\phi)\quad(\phi\in H_0^1(\Omega))
\]
for all $n\in \N$.

Then $u_n \rightharpoonup u\in  H_0^1(\Omega)$ and $a_n\grad u_n \rightharpoonup a\grad u \in L^2(\Omega)^N$, where $u$ is the solution of
\[
            \langle a^\hom \grad u , \grad \phi\rangle = f(\phi)\quad(\phi\in H_0^1(\Omega))
\]with $a^\hom \in \C^{N\times N}$ given by
\[
    a^\hom e_i\cdot e_j \coloneqq \int_{[0,1)^d} a(y)(e_i+\grad \chi_i(y))\cdot (e_j+\grad \chi_j(y))\d y,
\]
where $\chi_i \in H^1([0,1)^d)$ with $\chi_i\bot1$ and satisfying periodic boundary conditions and
\[
   \dive a (e_i + \grad \chi_i)=0\quad(i\in \{1,\ldots, N\}.
\]
\end{theorem}

When the weak operator topology is concerned, the next result is just a reformulation of the previous one. For this, we denote for a coefficient matrix $a$, the Lax--Milgram solution operator mapping from $H^{-1}(\Omega)$ to $H_0^1(\Omega)$ with the coefficient $a$, by $\mathcal{LM}(a) \in \mathcal{B}(H^{-1}(\Omega),H_0^1(\Omega))$. We obtain:

\begin{theorem} Let $\Omega\subseteq \R^N$ open, bounded. Then
\begin{align*}
 \mathcal{LM}(a_n) & \to \mathcal{LM}(a^\hom) \in  \mathcal{B}_\w(H^{-1}(\Omega),H_0^1(\Omega)). \\
 a_n\grad \mathcal{LM}(a_n) & \to a^\hom \grad \mathcal{LM}(a^\hom) \in \mathcal{B}_\w(H^{-1}(\Omega),L^2(\Omega)^N).
\end{align*}
\end{theorem}

\emph{Nonperiodic coefficients -- symmetric case.} The abstract definition of homogenisation problems for real and symmetric coefficients goes back to Spagnolo, see \cite{Spagnolo1967}. Note that originally the notion was used for the heat equation and that the `$G$' in the definition stands for `Green's operators'. Note that the following notion is a straightforward generalisation of the periodic case. To start off, let $\Omega\subseteq \R^N$ open and bounded. We introduce a particular class of coefficients
\[
   M_{\sym}(\alpha,\beta,\Omega)\coloneqq \{ a\in L^\infty(\Omega;\R^{N\times N}); a \text{ symmetric}, a(x)\geq \alpha, a(x)^{-1}\geq 1/\beta\, (\text{a.e. } x)\}.
\]
\begin{definition}[{{$G$-convergence, \cite[p. 476]{Spagnolo1976}}}] Let $(a_n)_n$ in $M_{\sym}(\alpha,\beta,\Omega)$ and $a\in M_{\sym}(\alpha,\beta,\Omega)$. Then $(a_n)_n$ is said to \emph{$G$-converge} to $a$, $a_n \stackrel{G}{\to} a$, if for all $f\in H^{-1}(\Omega)$ and $u_n\in H_0^1(\Omega)$ the following implication holds:
\[
     \langle a_n \grad u_n , \grad \phi\rangle = f(\phi)\quad(\phi\in H_0^1(\Omega))
\]implies $u_n\rightharpoonup u \in H_0^1(\Omega)$, where $u$ satisfies
\[
  \langle a \grad u , \grad \phi\rangle = f(\phi)\quad(\phi\in H_0^1(\Omega)).
\]
\end{definition}
We immediately read off that
\[
  (a_n)_n  \Gto a \iff \mathcal{LM}(a_n) \to \mathcal{LM}(a) \in  \mathcal{B}_\w(H^{-1}(\Omega),H_0^1(\Omega)).
\] Note that $G$-convergence of symmetric matrix coefficients induces a topology on $M_{\sym}(\alpha,\beta,\Omega)$, see \cite[Remark 2]{Spagnolo1976}. Call this topology $\tau_G$. Then we have the following result.
\begin{theorem}[{{\cite[Theorem 5]{Spagnolo1976}}}] $(M_{\sym}(\alpha,\beta,\Omega),\tau_G)$ is a metrisable, compact Hausdorff space.
\end{theorem}

This result can only be true, however, for the case of symmetric coefficients. this is confirmed by the next example.
\begin{example} Let $N=2$ and define $a\coloneqq 1 + \begin{pmatrix} 0 & -1 \\ 1 & 0\end{pmatrix}$. Then for all $\phi,\psi\in C_c^\infty(\Omega)$ we have
\begin{align*}
    \langle a \grad \phi,\grad\psi \rangle &=     \langle \grad \phi,\grad\psi \rangle - \langle \partial_2 \phi, \partial_1 \psi\rangle + \langle \partial_1 \phi,\partial_2\psi\rangle \\ & = \langle \grad \phi,\grad\psi \rangle + \langle  \phi,\partial_2 \partial_1 \psi \rangle+ \langle \partial_1 \phi,\partial_2\psi\rangle \\ & = \langle \grad \phi,\grad\psi \rangle - \langle  \partial_1\phi,\partial_2  \psi \rangle+ \langle \partial_1 \phi,\partial_2\psi\rangle \\ & = \langle \grad \phi,\grad\psi \rangle.
\end{align*}
Using that $C_c^\infty(\Omega)$ is dense in $H_0^1(\Omega)$, we deduce that
\[
     \langle a \grad \phi,\grad\psi \rangle = \langle \grad \phi,\grad\psi \rangle\quad (\phi,\psi\in H_0^1(\Omega)).
\]Hence, the variational formulation does not uniquely identify the coefficient, if the coefficient is allowed to be non-symmetric, as well. 
\end{example}
\emph{Nonperiodic coefficients -- nonsymmetric case.} The lack of separation can only be overcome, if also the fluxes are considered. For this, we set
\[
 M(\alpha,\beta,\Omega)\coloneqq \{ a\in L^\infty(\Omega)^{N\times N}; \Re a(x)\geq \alpha, \Re a(x)^{-1}\geq 1/\beta \ (\text{a.e. }x\in\Omega) \}.
\]
With this set of admissible coefficients, we introduce local $H$-convergence, invented by Tartar and Murat, where the $H$ stands for `homogenised'. 
\begin{definition}[{{$H$-convergence, see e.g.~\cite[Definition 6.4]{Tartar2009} or \cite{Murat1997}}}] Let $(a_n)_n$, $a$ in $M(\alpha,\beta,\Omega)$. Then $(a_n)_n$ \emph{$H$-converges to} $a$, $(a_n)_n \stackrel{H}{\to} a$, if for all $f\in H^{-1}(\Omega)$ and $u_n\in H_0^1(\Omega)$ the following implication holds:
\[
     \langle a_n \grad u_n , \grad \phi\rangle = f(\phi)\quad(\phi\in H_0^1(\Omega))
\]implies $u_n\rightharpoonup u \in H_0^1(\Omega)$ and $a_n\grad u_n\rightharpoonup a\grad u\in L^2(\Omega)^N$, where $u$ satisfies
\[
  \langle a \grad u , \grad \phi\rangle = f(\phi)\quad(\phi\in H_0^1(\Omega)).
  \]
\end{definition}
Again, one confirms immediately that $(a_n)_n\Hto a$, if and only if
\begin{align*}
 \mathcal{LM}(a_n) & \to \mathcal{LM}(a) \in  \mathcal{B}_\w(H^{-1}(\Omega),H_0^1(\Omega)) \text{ and }\\
 a_n\grad \mathcal{LM}(a_n) & \to a \grad \mathcal{LM}(a^\hom) \in \mathcal{B}_\w(H^{-1}(\Omega),L^2(\Omega)^N).
\end{align*}
Also, the convergence in the $H$-sense induces a topology, see \cite[below Definition 6.4]{Tartar2009}. We denote this topology on $M(\alpha,\beta,\Omega)$ by $\tau_H$. 
\begin{theorem}[{{Tartar, Murat, see e.g.~\cite[Theorem 6.5]{Tartar2009}}}] $(M(\alpha,\beta,\Omega),\tau_G)$ is a metrisable, compact Hausdorff space. Moreover, we have
\[
    (M_{\sym}(\alpha,\beta,\Omega),\tau_H)=(M_{\sym}(\alpha,\beta,\Omega),\tau_G).
\]
\end{theorem}

The above characterisations of the topology induced by $G$- and $H$-convergence can be slightly more refined. This, however, necessitates a closer look at $\mathcal{LM}(a)$; in fact $\mathcal{LM}(a)$ can be decomposed into three isomorphisms only one of which depending on $a$ and the others are mere implementations of the differential operators involved. The precise statement has been shown in \cite{TW14_FE} and reads as follows.

\begin{theorem}[{{\cite[Theorem 3.1]{TW14_FE} or \cite[Theorem 2.9]{W18_NHC}}}]\label{thm:TW} Let $H_0$ and $H_1$ be Hilbert spaces and let $C\colon \dom(C)\subseteq H_0\to H_1$ be a densely defined, closed linear operator with closed range; denote $\mathcal{C}\colon \dom(C)\cap \kar(C)^\bot\subseteq \kar(C)^\bot\to \rge(C), \phi\mapsto C\phi$ and endow the space $\dom(\mathcal{C})$ with its graph norm.

Let $a\in \mathcal{B}(H_1)$ satisfy
\[ 
    \Re \langle a \phi,\phi\rangle\geq \alpha\langle \phi,\phi\rangle\quad (\phi\in \dom(\mathcal{C}))
\]for some $\alpha>0$.

Then for all $f\in \dom(\mathcal{C})^*$ there exists a unique $u\in \dom(\mathcal{C})$ such that
\[
     \langle a Cu , C\phi\rangle = f(\phi)\quad (\phi\in \dom(\mathcal{C})).
\]
More precisely, we have
\[
     u =   \mathcal{C}^{-1} ( \iota^* a \iota)^{-1}  (\mathcal{C}^\diamond)^{-1}f,
\]
where $\mathcal{C}^\diamond \colon \rge(C) \to \dom(\mathcal{C})^*, q\mapsto (\phi\mapsto \langle q, C\phi\rangle_{H_1})$ and $\iota\colon \rge(C)\hookrightarrow H_1$ is the canonical embedding.
\end{theorem}
\begin{remark}
In Theorem \ref{thm:TW}, if it was not for the solution formula for $u$, the statement of the theorem would have been covered by the classical Lax--Milgram lemma. The solution formula, however, is the decisive point needed for refining the characterisations of $G$- and $H$-convergence above.
\end{remark}
\begin{example} (a) An elementary example for Theorem \ref{thm:TW} would be $C=\grad$ with $\dom(C)=H_0^1(\Omega)$, $H_0=L^2(\Omega)$ and $H_1=L^2(\Omega)^N$ and any strictly positive definite operator $a\in \mathcal{B}(L^2(\Omega)^N)$. Note that in this case $\mathcal{C}=\grad:H_0^1(\Omega)\to \grad[H_0^1(\Omega)]$. Recall that $\grad[H_0^1(\Omega)]$ is closed due to Poincar\'e's inequality. Moreover, we have $\mathcal{C}^\diamond = \dive \colon \grad[H_0^1(\Omega)]\to H^{-1}(\Omega)$. We note that this is the standard distributional divergence operator applied to general $L^2(\Omega)^N$ vector fields with the nullspace of the standard $\dive$ being removed.

(b) A more involved example is the $\curl$-operator with homogeneous `electric' boundary conditions as a particular variant of $C$. In this case $H_0=H_1=L^2(\Omega)^3$.
\end{example}
Using that the gradient and the divergence realised as $\mathcal{C}$ and $\mathcal{C}^\diamond$, respectively, are isomorphism, we straightforwardly obtain the desired more refined characterisations.
\begin{theorem}[{{\cite[Remark 4.11]{W18_NHC} and \cite[Theorem 1.2]{W16_Gcon}}}]\label{thm:GH} Let $\Omega\subseteq \R^N$ be open and bounded and denote $\iota_{g_0} \colon g_0(\Omega) \to L^2(\Omega)^N$ the canonical embedding, where $g_0(\Omega)\coloneqq \grad[H_0^1(\Omega)]$.

(a) Let $(a_n)_n$ and $a$ in $M_{\sym}(\alpha,\beta,\Omega)$. Then
\[
  (a_n)_n \Gto a \iff  \left((\iota_{g_0}^* a_n \iota_{g_0})^{-1}\right)_n \to (\iota_{g_0}^* a \iota_{g_0})^{-1} \in \mathcal{B}_\w(g_0(\Omega)).
\]
(b)  Let $(a_n)_n$ and $a$ in $M(\alpha,\beta,\Omega)$. Then
\[
  (a_n)_n \Hto a \iff  \begin{cases} \left((\iota_{g_0}^* a_n \iota_{g_0})^{-1}\right)_n \to (\iota_{g_0}^* a \iota_{g_0})^{-1} \in \mathcal{B}_\w(g_0(\Omega)) & \text{and}\\
                   a_n\iota_{g_0} \left((\iota_{g_0}^* a_n \iota_{g_0})^{-1}\right)_n \to a\iota_{g_0} (\iota_{g_0}^* a \iota_{g_0})^{-1} \in \mathcal{B}_\w(g_0(\Omega),L^2(\Omega)^N). &
  \end{cases}
\]
\end{theorem}
\begin{remark}
A related characterisation can be found in \cite[Example 6.7]{EGW17_D2N}. Note that in \cite{EGW17_D2N} convergence of the Dirichlet-to-Neumann operator under $G$-convergence has been addressed.
\end{remark}

\section{Time-independent problems -- nonlocal coefficients}

\emph{An example for a homogenisation problem with nonlocal coefficients.} Similar to the local coefficient case, the story of homogenisation problems starts with a particular model problem. For this let $k\in  L^\infty(\R^3)$ such that $k(x+\ell)=k(x)$ for a.e. $x\in \R^3$ and all $\ell\in \Z^3$. Then we let $k_n \coloneqq k(n\cdot)$ for all $n\in \N$ and consider
\[
      k_n * f\coloneqq \left( x \mapsto \int_\Omega k_n(x-y)f(y)\d y\right).
\]
Let $\Omega\subseteq \R^3$ be bounded, simply connected, weak Lipschitz domain with connected complement.
Assume that $\|k\|_{L^\infty(\R^3)}<\lambda(\Omega)$. This then implies that $\sup_{n\in\N}\|k_n\|_{L^1(\Omega)}<1$ and hence $\sup_{n\in\N}\|k_n*\|_{\mathcal{B}(L^2(\Omega))}<1$, by Young's inequality. Consequently, $\Re (1-k_n*)\geq \alpha$ for some $\alpha>0$ uniformly in $n$. The model for nonlocal stationary diffusion is now to find $u_n\in H_0^1(\Omega)$ such that for $f\in H^{-1}(\Omega)$ we have
\[
   -\dive(1-k_n*)\grad u_n = f.
\]
As before, we consider the limit $n\to\infty$ and obtain $u_n \rightharpoonup u \in H_0^1(\Omega)$, where $u$ satisfies
\[
    -\dive(1-\mathfrak{M}(k)\chi_{\Omega}*)\grad u = f.,
\]
where $\chi_{\Omega}*$ denotes convolution with the characteristic function of $\Omega$.
A proof of this fact will be given below. Obviously, the just mentioned result is not covered by local coeffcients. In the following, we will see that in order to deduce the said limit behaviour, one needs to look beyond classical divergence form problems. A major step will be the next observation.

\emph{Helmholtz decompositions.} In order to deal with nonlocal coefficients also on a more general level, we decompose the underlying $L^2$-space. For any open set $\Omega\subseteq \R^3$, we recall and denote
\begin{align*}
    g_0(\Omega) & = \{\grad u; u\in H_0^1(\Omega)\} \\
    g(\Omega) & = \{\grad u; u\in H^1(\Omega)\} \\
    c(\Omega) & =  \{\curl u; u\in H(\curl)\} \\
        c_0(\Omega) & =  \{\curl u; u\in H_0(\curl)\} 
\end{align*}
\begin{theorem}[\cite{Picard1984,Picard1982,Picard1983,Picard1990}]\label{thm:HD}Let $\Omega\subseteq \R^3$ be a bounded, open, weak Lipschitz domain. Then we have
\[
    L^2(\Omega)^3 = g_0(\Omega)\oplus c(\Omega)\oplus \mathcal{H}_\textrm{D}(\Omega) = g(\Omega)\oplus c_0(\Omega)\oplus \mathcal{H}_\textrm{N}(\Omega),
\]
where \begin{align*}
 \mathcal{H}_\textrm{D}(\Omega)& =\{ u\in H_0(\curl); \dive u = 0, \curl u =0\}, \\
 \mathcal{H}_\textrm{N}(\Omega)& =\{ u\in H_0(\dive); \dive u = 0, \curl u =0\} 
 \end{align*}
 are the harmonic Dirichlet and Neumann fields, respectively. Note that the dimension of $\mathcal{H}_{\textrm{N}}(\Omega)$ equals the number of pairwise non-homotopic closed curves not homotopic to a single point; the dimension of $\mathcal{H}_{\textrm{D}}(\Omega)$ counts the number of connected components of the complement of $\Omega$ ignoring the component of the exterior domain surrounding $\Omega$.
\end{theorem}
\begin{corollary}\label{cor:toptriv}
Let $\Omega \subseteq \R^3$ bounded, open, weak Lipschitz domain. If $\Omega$ is simply connected, then
\[
   L^2(\Omega)^3 = g(\Omega)\oplus c_0(\Omega).
\]
If $\Omega$ has connected complement, then
\[
   L^2(\Omega)^3 = g_0(\Omega)\oplus c(\Omega).
\]
\end{corollary}

\begin{definition} We call  $\Omega$ \emph{topologically trivial}, if $\dim(\mathcal{H}_\textrm{D}(\Omega))=\dim(\mathcal{H}_\textrm{N}(\Omega))=0$; $\Omega\subseteq \R^3$ is called \emph{standard domain}, if $\Omega$ is open, bounded, weak Lipschitz and topologically trivial with the segment property.
\end{definition}

\emph{Block operator matrix representation.} Let $\Omega$ be a standard domain. A key to most conveniently describe actions of operators in $L^2(\Omega)^3$ is by using the decomposition result in Corollary \ref{cor:toptriv}. We define $\iota_{g_0}$ to be the canonical embedding from $g_0(\Omega)$ into $L^2(\Omega)^3$; we also define $\iota_c$, etc.~similarly. For any $a\in \mathcal{B}(L^2(\Omega)^3)$ we put
\[
    \begin{pmatrix}
       a_{00} & a_{01} \\ a_{10} & a_{11}
    \end{pmatrix} \coloneqq     \begin{pmatrix}
       \iota_{g_0}^* a\iota_{g_0} & \iota_{g_0}^*a\iota_{c}\\ \iota_{c}^*a\iota_{g_0} & \iota_{c}^*a\iota_{c}
    \end{pmatrix}.
\]
Analogously, we set,
\[
    \begin{pmatrix}
       \hat{a}_{00} & \hat{a}_{01} \\ \hat{a}_{10} & \hat{a}_{11}
    \end{pmatrix} \coloneqq     \begin{pmatrix}
       \iota_{g}^* a\iota_{g} & \iota_{g}^*a\iota_{c_0}\\ \iota_{c_0}^*a\iota_{g} & \iota_{c_0}^*a\iota_{c_0}
    \end{pmatrix}.
\]

\emph{Definition of nonlocal $H$-convergence.} Equipped with the latter notion and notation, we are able to define nonlocal $H$-convergence. In fact, there is an abundance of choice which decomposition to use. It can be shown that the two possibilities provided lead to different topologies, that is, nonlocal $H$-convergence `sees' the boundary conditions. This is one central property, which local $H$- or $G$-convergence do not share with nonlocal $H$-convergence. Let $\Omega$ be a standard domain. The set of admissible coefficients is given as follows
\begin{align*}
    \mathcal{M}(\alpha,\beta,\Omega) & \coloneqq \{ a\in \mathcal{B}(L^2(\Omega)^3); a \text{ invertible},\\ &\hspace*{3cm} \Re a_{00} \geq \alpha, \Re a_{00}^{-1}\geq 1/\beta, \Re (a^{-1})_{11}\geq 1/\beta, \Re (a^{-1})_{11}^{-1}\geq \alpha \}\\
        \hat{\mathcal{M}}(\alpha,\beta,\Omega)& \coloneqq \{ a\in \mathcal{B}(L^2(\Omega)^3); a \text{ invertible},\\ &\hspace*{3cm} \Re \hat{a}_{00} \geq \alpha, \Re \hat{a}_{00}^{-1}\geq 1/\beta, \Re (\hat{a^{-1}})_{11}\geq 1/\beta, \Re (\hat{a^{-1}})_{11}^{-1}\geq \alpha \}.
\end{align*} We note here that $M(\alpha,\beta,\Omega)\subseteq   \mathcal{M}(\alpha,\beta,\Omega)\cap   \hat{\mathcal{M}}(\alpha,\beta,\Omega)$.
  In order to fit to a more general perspective, see below, we shall define nonlocal $H$-convergence as follows. In contrast to the introduction of $H$-convergence and $G$-convergence, we define the topology first, see \cite{W18_NHC}.
\begin{definition} (a) The \emph{topology of nonlocal $H$-convergence w.r.t.~$(g_0(\Omega),c(\Omega))$} is the initial topology $\tau_{\nlH}$ on $\mathcal{M}(\alpha,\beta,\Omega)$ such that
\begin{align*}
   a& \mapsto a_{00}^{-1} \in \mathcal{B}_\w (g_0(\Omega)) \\
   a& \mapsto a_{10}a_{00}^{-1} \in \mathcal{B}_\w (g_0(\Omega),c(\Omega)) \\
   a& \mapsto a_{00}^{-1} a_{01} \in \mathcal{B}_\w (c(\Omega),g_0(\Omega)) \\
   a& \mapsto a_{11}-a_{10}a_{00}^{-1} a_{01} \in \mathcal{B}_\w (c(\Omega)) \\
\end{align*}
are continuous.

(b) The \emph{topology of nonlocal $H$-convergence w.r.t.~$(g(\Omega),c_0(\Omega))$ }is the initial topology $\hat{\tau}_{\nlH}$ on $\hat{\mathcal{M}}(\alpha,\beta,\Omega)$ such that
\begin{align*}
   a& \mapsto \hat{a}_{00}^{-1} \in \mathcal{B}_\w (g(\Omega)) \\
   a& \mapsto \hat{a}_{10}\hat{a}_{00}^{-1} \in \mathcal{B}_\w (g(\Omega),c_0(\Omega)) \\
   a& \mapsto \hat{a}_{00}^{-1} \hat{a}_{01} \in \mathcal{B}_\w (c_0(\Omega),g(\Omega)) \\
   a& \mapsto \hat{a}_{11}-\hat{a}_{10}\hat{a}_{00}^{-1} \hat{a}_{01} \in \mathcal{B}_\w (c_0(\Omega)) \\
\end{align*}
are continuous.
\end{definition}
\begin{remark}
From the definition of the topology it is unclear, whether the summing of operators is continuous under this topology. In fact, already one-dimensional (and local) analogues  of the topology just introduced shows that addition is not continuous. However, the multiplication by scalars is continuous.
\end{remark}

Using that the unit ball of the space of bounded linear operators acting on a separable Hilbert space is compact, Hausdorffian and metrisable under the weak operator topology, we immediately obtain the following result.
\begin{theorem}[{{\cite[Theorem 5.5]{W18_NHC}}}] Let $\mathcal{B}\subseteq \mathcal{M}(\alpha,\beta,\Omega)$ be bounded. Then $(\overline{\mathcal{B}}^{\tau_{\nlH}},\tau_{\nlH})$ is compact and metrisable.
\end{theorem}
Of course, a similar result holds for the other stated variant of nonlocal $H$-convergence.
For completeness, we just recall that the mentioned variants of nonlocal $H$-convergence cannot be compared.
\begin{proposition}[{{\cite[Example 5.3]{W18_NHC}}}]
Consider
\begin{align*}
   I \colon  \left( \mathcal{M}(\alpha,\beta,\Omega)\cap \hat{\mathcal{M}}(\alpha,\beta,\Omega), \tau_{\nlH} \right) & \to  \left( \mathcal{M}(\alpha,\beta,\Omega)\cap \hat{\mathcal{M}}(\alpha,\beta,\Omega), \hat{\tau}_{\nlH} \right) \\
   a& \mapsto a.
\end{align*}
Then neither $I$ nor $I^{-1}$ is continuous.
\end{proposition}

\emph{Characterisation with solutions to PDEs.} Although the previous paragraph is formally at the heart of this survey in as much as it heavily involves the weak operator topology, the link to the traditional notions of convergence is still missing. A fundamental step towards this aim is the following characterisation. Let $\Omega$ be a standard domain. We focus on $(\mathcal{M}(\alpha,\beta,\Omega),\tau_{\nlH})$; the result for $(\hat{\mathcal{M}}(\alpha,\beta,\Omega),\hat{\tau}_{\nlH})$ is similar. We need the space $\tilde{H}(\curl)\coloneqq \{ q\in H(\curl); q\in c_0(\Omega)\}$ and $\tilde{H}^{-1}(\curl)\coloneqq \tilde{H}(\curl)^*$.
\begin{theorem}[{{\cite[Theorem 4.1]{W18_NHC}}}] \label{thm:chPDE} Let $(a_n)_n$ and $a$ in  $\mathcal{M}(\alpha,\beta,\Omega)$. Then $(a_n)_n\to a$ in $\tau_{\nlH}$ if and only if the following statement holds:  for all $f\in H^{-1}(\Omega)$ and $g\in \tilde{H}^{-1}(\curl)$ let $u_n\in H_0^1(\Omega)$ and $v_n\in \tilde{H}(\curl)$ satisfy 
\begin{align*}
  \langle a_n \grad u_n ,\grad \phi\rangle = f(\phi) & \quad (\phi\in H_0^1(\Omega)) \\
  \langle a_n^{-1} \curl v_n ,\curl \psi\rangle =g(\psi) & \quad (\psi\in \tilde{H}(\Omega)).
\end{align*}
Then $u_n \rightharpoonup u \in H_0^1(\Omega)$ and $v_n\rightharpoonup v\in \tilde{H}(\curl)$ and $a_n\grad u_n \rightharpoonup a\grad u$ and $a_n^{-1}\curl v_n \rightharpoonup a^{-1}\curl v$, where $u$ and $v$ satisfy
\begin{align*}
  \langle a \grad u ,\grad \phi\rangle = f(\phi) & \quad (\phi\in H_0^1(\Omega)) \\
  \langle a^{-1} \curl v ,\curl \psi\rangle =g(\psi) & \quad (\psi\in \tilde{H}(\curl)).
\end{align*}
\end{theorem}
This characterisation provides the desired connections to PDEs. When it comes to practical applications, it is helpful to have the following result at hand.
\begin{theorem}[{{\cite[Theorem 6.2 (also cf.~Theorem 6.5)]{W18_NHC}}}]\label{thm:DCLch} Let $(a_n)_n$ and $a$ in  $\mathcal{M}(\alpha,\beta,\Omega)$. Then the following conditions are equivalent:

(a)  $(a_n)_n\to a$ in $\tau_{\nlH}$

(b) for all weakly convergent $(q_n)_n$ in $L^2(\Omega)^3$, $q=\w\text{-}\lim_{n\to\infty} q_n$, and $\kappa\colon \N\to \N$ strictly monotone the following implication holds:
Assume that $(\dive a_{\kappa(n)}q_n)_n$ is relatively compact in $H^{-1}(\Omega)$ and $(\curl q_n)_n$ is relatively compact in $\tilde{H}^{-1}(\curl)$. Then $a_{\kappa(n)}q_n \rightharpoonup aq$ as $n\to\infty$.
\end{theorem}

The latter characterisation can be used to verify the convergence statement, where we treated the nonlocal coefficient above. For this we need the following form of a global $\dive$-$\curl$-lemma; see also \cite{Pauly2017}.
\begin{theorem}[{{\cite[Theorem 2.6]{W17_DCL}}}]\label{thm:DCL} Let $\Omega\subseteq \R^3$ be an open, bounded weak Lipschitz domain. Let $(q_n)_n, (r_n)_n$ in $L^2(\Omega)^3$ weakly convergent. Assume that $(\dive r_n)_n$ is relatively compact in $H^{-1}(\Omega)$ and that $(\curl q_n)_n$ is weakly convergent in $\tilde{H}^{-1}(\curl)$.

 Then $\langle q_n,r_n\rangle \to \langle \w\text{-}\lim q_n,\w\text{-}\lim r_n\rangle$ as $n\to\infty$. 

\end{theorem}

\begin{example} Let $\Omega\subseteq \R^3$ be a standard domain. Then $(1-k_n*)\to (1-\mathfrak{M}(k)\chi_{\Omega}*)$ in $\tau_{\nlH}$, where $k \in L^\infty(\R^3)$ is $[0,1)^3$-periodic and $k_n\coloneqq k(n\cdot)$ and $\mathfrak{M}(k)=\int_{[0,1)^3}k$.  We use Theorem \ref{thm:DCLch} in order to show the desired convergence statement. For this let $(q_n)_n$ in $L^2(\Omega)^3$ be weakly convergent to some $q\in L^2(\Omega)^3$ and $\kappa\colon \N\to \N$ strictly monotone. Assume that $(\dive(1-k_n*)q_n)_n$ is relatively compact in $H^{-1}(\Omega)$ and that $(\curl q_n)_n$ is relatively compact in $\tilde{H}^{-1}(\curl)$. We need to show that $(1-k_{\kappa(n)}*)q_n \rightharpoonup (1-\mathfrak{M}(k)\chi_{\Omega}*)q.$ For convenience, we drop $\kappa(n)$ and write $n$ instead. Note that $((1-k_{n}*)q_n)_n$ is uniformly bounded in $L^2(\Omega)^3$, so we need to check weak convergence only on a dense subset. For this, let $\phi\in C_c^\infty(\Omega)^3$ and consider
\[
     \langle (1-k_{n}*)q_n, \phi\rangle =      \langle q_n, (1-k_{n}*)^*\phi\rangle =   \langle q_n, \phi-(k_{n}*)^*\phi\rangle.
\]
It is elementary to see (use periodicity of $k$; see also the introduction) that $r_n \coloneqq (k_{n}*)^*\phi \rightharpoonup \mathfrak{M}(k^*)(\chi_{\Omega}*)^*\phi$. By assumption, we have $(\curl q_n)_n$ is relatively compact in $\tilde{H}^{-1}(\curl)$. Next $\dive r_n = (k_{n}*)^*\dive \phi \in L^2(\Omega)^3$. As $\Omega$ is bounded, $L^2(\Omega)^3$ embeds compactly into $H^{-1}(\Omega)$. Hence, $(\dive r_n)_n$ is relatively compact in $H^{-1}(\Omega)$. Thus, by Theorem \ref{thm:DCL}, we deduce
\[
     \langle q_n, \phi-(k_{n}*)^*\phi\rangle\to   \langle q, \phi-\mathfrak{M}(k^*)(\chi_{\Omega}*)^*\phi\rangle = \langle (1-\mathfrak{M}(k)\chi_{\Omega}*)q,\phi\rangle,
\]
which eventually yields the assertion.
\end{example}

\emph{Connections to the local topologies.} The characterisation in Theorem \ref{thm:chPDE} already shows that a sequence in $M(\alpha,\beta,\Omega)$ or $M_{\sym}(\alpha,\beta,\Omega)$, which is nonlocally $H$-convergent, yields an $H$- or $G$-convergent sequence. The limits in $\tau_{\nlH}$ and w.r.t.~$H$- or $G$-convergence coincide. In fact, even more is true:
\begin{theorem}[{{\cite[Theorem 5.11 and Remark 5.12]{W18_NHC}}}]\label{thm:locnonloc} Let $\Omega\subseteq \R^3$ be a standard domain. Then 
\[
   (M(\alpha,\beta,\Omega),\tau_{\nlH})=   (M(\alpha,\beta,\Omega),\hat{\tau}_{\nlH})=   (M(\alpha,\beta,\Omega),\tau_{H})
\]
and
\[
   (M_{\sym}(\alpha,\beta,\Omega),\tau_{\nlH})=   (M_{\sym}(\alpha,\beta,\Omega),\hat{\tau}_{\nlH})=   (M_{\sym}(\alpha,\beta,\Omega),\tau_{G}).
\]
\end{theorem}

\emph{Generalisations.} In this section we have focussed on the three-dimensional case and problems of `classical' divergence form. The whole concept of nonlocal $H$-convergence, however, has been developed for closed Hilbert complexes. More precisely, the operators $\dive$, $\grad$, and $\curl$ are suitably replaced by closed linear operators $A_2, A_1, A_0$ with closed ranges and the property that $\rge(A_0)=\kar(A_1)$ and $\rge(A_1)=\kar(A_2)$. Thus, the concepts developed above also work for other (also mixed type) boundary conditions. Moreover, other equations can be considered as well. The developed concepts naturally work for elasticity of the biharmonic operator, using the Pauly--Zulehner complex, see \cite{Pauly2019}.

\section{Dynamic Problems -- Preliminaries}\label{s:set}

\emph{Time-derivative and exponentially weighted $L^2$-spaces.} We briefly introduce the concept of evolutionary equations and the operator-theoretic notions accompanied by these. For all the statements in this section we refer to \cite{A11,KPSTW14_OD,W16_H} for a more detailed exposition. The original paper is \cite{PicPhy}. Throughout, let $\nu\in\mathbb{R}$  and $H$ be a Hilbert space. Denote by $L_\nu^2(\mathbb{R};H)$ the Hilbert space of (equivalence classes of) Bochner measurable functions $f\colon\mathbb{R}\to H$ such that
\[
   \|f\|_{L_\nu^2}^2=\langle f,f\rangle_{L_\nu^2}=\int_\mathbb{R} \|f(t)\|_H^2 \exp(-2t\nu)\dd t
\]
is finite. Denoting by $f'$ the distributional derivative of $f\in L_{\textnormal{loc}}^1(\mathbb{R};H)$, we define
\[
   \partial_{t,\nu} \colon H_\nu^1(\mathbb{R};H) \subseteq L_\nu^2(\mathbb{R}:H)\to L_\nu^2(\mathbb{R};H), f\mapsto f',
\]
where $H_\nu^1(\mathbb{R};H)$ denotes the Sobolev space of $L_\nu^2(\mathbb{R};H)$-functions with $f'\in L_\nu^2(\mathbb{R};H)$. 

Next, we introduce the \emph{Fourier--Laplace transformation}. Define $\mathcal{L}_\nu$ by
\[
   \mathcal{L}_\nu \phi (\xi)\coloneqq \frac{1}{\sqrt{2\pi}}\int_\mathbb{R} \ee^{-\ii\xi t-\nu t}\phi(t)\dd t\quad(\xi\in\R),
\]
where $\phi\colon \mathbb{R}\to H$ is continuous and compactly supported. By a variant of Plancherel's theorem, $\mathcal{L}_\nu$ admits a unitary extension as an operator from $L_\nu^2(\mathbb{R};H)$ to $L^2(\mathbb{R};H)$. We shall re-utilise $\mathcal{L}_\nu$ for this extension. 

The Fourier--Laplace transformation yields a spectral representation of $\partial_t$. For this to make precise, we denote by
\begin{align*}
   \mm \colon \dom(\mm)\subseteq L^2(\mathbb{R};H) &\to L^2(\mathbb{R};H),
    \\ f&\mapsto (\xi \mapsto \xi f(\xi))
\end{align*}
the multiplication-by-the-argument operator, where $\dom(\mm)=\{f\in L^2(\mathbb{R};H); (\xi \mapsto \xi f(\xi))\in L^2(\mathbb{R};H)\}$. By a slight abuse of notation for the operator $\lambda I_H$, where $I_H$ is the identity in $H$,  we shall always just write $\lambda$. The explicit spectral theorem for $\partial_t$ now reads as follows.

\begin{theorem}[see e.g.~{{\cite[Corollary 2.5]{KPSTW14_OD}}}]\label{thm:spr}
  For all $\nu\in\mathbb{R}$, we have
  \[
     \partial_{t,\nu} = \mathcal{L}_\nu^*(\ii \mm + \nu)\mathcal{L}_\nu.
  \]
\end{theorem}

\emph{Material law operators and functions of $\partial_t$.}
Being unitarily equivalent to a normal operator, the operator $\partial_t$ is normal itself. More importantly, Theorem \ref{thm:spr} provides a functional calculus for $\partial_t$. We restrict ourselves to a class of holomorphic functions only. In fact, in order to obtain a so-called \emph{causal} operator, this restriction is necessary, see e.g. \cite[Corollary 1.2.5]{W16_H} or \cite{Weiss1991} for accessible proofs of a result due to \cite{Foures1955}. The space of analytic mappings from an open set $E\subseteq \mathbb{C}$ to some Banach space $X$ is denoted by $\mathcal{H}(E;X)$.

\begin{definition}
  Let $M\in \mathcal{H}(\mathbb{C}_{\Re>\nu};
  \mathcal{B}(K,H))$, $H,K$ Hilbert spaces. For all $\mu>\nu$, we define
  \[
     M(\partial_{t,\mu})\coloneqq \mathcal{L}_\mu^* M(\ii \mm+\mu)\mathcal{L}_\mu,
  \]
  where we endow $M(\ii \mm +\mu)$ with its maximal domain. 
\end{definition}

\begin{remark}\label{rem:ind}
  Note that $M(\partial_{t,\mu})$ is realised as a densely defined operator acting from $L_\mu^2(\mathbb{R};K)$ to $L_\mu^2(\mathbb{R};H)$ for all $\mu>\nu$. By \cite[Lemma 3.6]{Trostorff2013}, for $f\in \dom(M(\partial_{t,\mu}))\cap \dom(M(\partial_{t,\eta}))$ we have
  \[
     M(\partial_{t,\mu})f= M(\partial_{t,\nu})f.
  \] For this reason, we shall also employ the custom to dispense with mentioning $\mu$ or $\nu$.
\end{remark}

\emph{An abstract class of PDEs.} Next, we define the notion of evolutionary equations. We emphasise that the term `evolutionary' is used in order to distinguish from explicit Cauchy problems, which are commonly summarised by `evolution equations' and form a proper subclass of evolutionary equations. Let $H$ be a Hilbert space.
\begin{definition}
  Let $A\colon \dom(A)\subseteq H \to H$ be densely defined and closed, $M\in \mathcal{H}(\mathbb{C}_{\Re>\nu};\mathcal{B}(H))$. For $U,F\in L_\mu^2(\mathbb{R};H)$, $\mu>\nu$. An equation of the form
  \[
     (M(\partial_{t,\mu})+A)U=F
  \]
  is called \emph{evolutionary equation}. An evolutionary equation (or $M(\partial_{t,\mu})+A$)  is called \emph{well-posed}, if
  \[
     {S}\colon \mathbb{C}_{\Re>\nu} \ni z\mapsto (M(z)+A)^{-1}\in \mathcal{B}(H)
  \]
  is well-defined and we have
  \[
     {S}\in \mathcal{H}^\infty(\C_{\Re>\nu};\mathcal{B}(H))\coloneqq \{ T\in \mathcal{H}(\mathbb{C}_{\Re>\nu};\mathcal{B}(H)); T \text{ is bounded}\}.
  \]
\end{definition}
\begin{remark}
As a matter of jargon, due to Remark \ref{rem:ind}, we shall also say $M(\partial_{t})+A$ is well-posed in $\mathcal{B}(L^2_\mu(\R;H))$, by which we mean that $M\in \mathcal{H}(\mathbb{C}_{\Re>\nu};\mathcal{B}(H))$ for some $\nu<\mu$ and that $ \mathbb{C}_{\Re>\nu} \ni z\mapsto (M(z)+A)^{-1}\in \mathcal{B}(H)$ is well-defined and bounded.
\end{remark}

By the Fourier--Laplace transformation it is easy to see that for well-posed $M(\partial_{t,\mu})+A$, the operator $S(\partial_{t,\mu})=(M(\cdot)+A)^{-1}(\partial_{t,\mu})$ is a bounded linear operator in $L_\mu^2(\mathbb{R};H)$ for all $\mu>\nu$. 

\begin{remark} Let $A\colon \dom(A)\subseteq H\to H$ be a generator of a $C_0$-semigroup. By standard $C_0$-semigroup theory, see e.g.~\cite{Engel2000}, there exists $\omega_0\in\mathbb{R}$ such that $\mathbb{C}_{\Re>\omega_0}\subseteq \rho(A)$. Moreover, we have
\[
   \mathbb{C}_{\Re>\omega_0} \ni z \mapsto (z-A)^{-1} 
\]
is bounded and analytic. Hence, $(\partial_t + A)$ is well-posed. 
\end{remark}

The standard case of evolutionary equations has been introduced in \cite{PicPhy}.

\begin{theorem}[{{\cite[Solution Theory]{PicPhy}}}]\label{thm:st} Let $N\in \mathcal{H}^\infty(\mathbb{C}_{\Re>\nu};L(H))$, $A\colon \dom(A)\subseteq H\to H$ skew-self-adjoint. Assume there exists $c>0$ such that for all $z\in \mathbb{C}_{\Re>\nu}$
\[
   \Re z N(z) = \frac{1}{2}\left( z N(z)+z^* N(z)^*\right)\geq c.
\]
Then, for all $\mu>\nu$,  $M(\partial_{t,\mu})+A$ is well-posed, where $M(z)\coloneqq zN(z)$, $z\in\mathbb{C}_{\Re>\nu}$ and
\[
   S\coloneqq \overline{M(\partial_{t,\mu})+A}^{-1}=(M(\cdot)+A)^{-1}(\partial_{t,\mu}).
\] The \emph{solution operator} S is time-shift invariant and leaves functions supoorted on $[0,\infty)$ invariant, i.e., $S$ is \emph{causal}.
\end{theorem}
Next, we shall present several (standard) examples. We shall also refer to \cite[Section 2.1 Guiding Examples]{PTW15_WP_P} for a more detailed account on the equations to follow. In the lines to come as well as in the next sections, the time-derivative operator $\partial_t$ acts on the first variable only; we shall assume always implicit zero initial conditions at $-\infty$, the spatial operators $\dive, \grad$ and $\curl$ only act on the variables of $\Omega$.

\emph{Heat conduction.} Let $\Omega\subseteq \R^N$ be open. The equations for heat conduction are given by the heat balance law, which says
\[
    \partial_t \theta +\dive q = Q,
\]
accompanied by Fourier's law
\[
    q = - a\grad_0 \theta,
\]
where $\theta\colon \R\times \Omega\to \R$ is the unknown heat, $q \colon \R\times \Omega\to \R^N$ is the unknown heat flux, and $Q$ is the given heat source and $a\in \mathcal{B}(L^2(\Omega)^N)$ is the given heat conductivity. We shall assume that $\theta$ satisfies homogeneous Dirichlet boundary conditions; we stress this boundary condition to holds by writing $\grad_0$ instead of $\grad$. Note that the domain of definition of $\grad_0$ is $H_0^1(\Omega)$.
Assuming that $\Re a\geq \alpha>0$, we can apply the well-posedness Theorem \ref{thm:st} with $A=\begin{pmatrix} 0 & \dive \\ \grad_0 & 0 \end{pmatrix}$ and $M(z)=z \begin{pmatrix} 1 & 0 \\ 0 & 0 \end{pmatrix}+ \begin{pmatrix} 0 & 0 \\ 0 & a^{-1} \end{pmatrix}$. The boundary conditions for $\theta$ imply that $A$ is skew-self-adjoint.  Thus,
\[
    \partial_t \begin{pmatrix} 1 & 0 \\ 0 & 0 \end{pmatrix}+ \begin{pmatrix} 0 & 0 \\ 0 & a^{-1} \end{pmatrix} +\begin{pmatrix} 0 & \dive \\ \grad_0 & 0 \end{pmatrix}
\]is well-posed in $L^2_\mu(\R;L^2(\Omega)^{N+1})$ for all $\mu>0$.

\emph{Wave equation.} Let $\Omega\subseteq \R^N$ open. The scalar elastic or acoustic wave equation is given by a balance of momentum equation
\[
    \partial_t^2 u - \dive \sigma = F,
\]
accompanied with
\[
     \sigma = a\grad_0 u.
\]
Here, $u\colon \R\times \Omega \to \R$ and $\sigma\colon \R\times \Omega \to \R^N$, the displacement and the stress, respectively, are the unknowns and $F$ is the given elastic force as well as $a=a^*\geq \alpha>0$ is the given elasticity tensor. Again, we shall assume homogeneous Dirichlet boundary condition for $u$.  Substituting $v=\partial_t u$, we obtain with $A=\begin{pmatrix} 0 & \dive \\ \grad_0 & 0 \end{pmatrix}$ and  $M(z)=z \begin{pmatrix} 1 & 0 \\ 0 & a^{-1} \end{pmatrix}$ and the well-posedness Theorem \ref{thm:st} that
\[
    \partial_t \begin{pmatrix} 1 & 0 \\ 0 & a^{-1} \end{pmatrix} +\begin{pmatrix} 0 & \dive \\ \grad_0 & 0 \end{pmatrix}
\]is well-posed in $L^2_\mu(\R;L^2(\Omega)^{N+1})$ for all $\mu>0$.

\emph{Maxwell's equations.} Let $\Omega\subseteq \R^3$ be open. Maxwell's equations in matter are given by the following two equations:
\begin{align*}
  \partial_t \varepsilon E + \sigma E -\curl H & = -J \\
    \partial_t \mu H + \curl_0 E & = 0,
\end{align*}
where $E$ satisfies the homogeneous electric boundary condition of vanishing tangential component at the boundary (stressed by writing $\curl_0$). The unknown in Maxwell's equations is the electro-magnetic field $(E,H)$. The given material-dependent quantities are $\epsilon=\epsilon^*\geq \alpha,\mu=\mu^*\geq \alpha ,\sigma\in \mathcal{B}(L^2(\Omega)^3)$  the dielelctricity, the magnetic permittivity and the electric conductivity; the right-hand side $-J$ is a given forcing term due to external currents. With the setting $M(z) = z \begin{pmatrix} \varepsilon & 0 \\ 0 & \mu \end{pmatrix} + \begin{pmatrix} \sigma & 0 \\ 0 & 0 \end{pmatrix}$ and $A=\begin{pmatrix} 0 & -\curl \\ \curl_0 & 0 \end{pmatrix}$ (note that the electric boundary condition for the lower left $\curl$-operator makes $A$ skew-selfadjoint), we obtain
\[
    \partial_t \begin{pmatrix} \varepsilon & 0 \\ 0 & \mu \end{pmatrix} + \begin{pmatrix} \sigma & 0 \\ 0 & 0 \end{pmatrix} + \begin{pmatrix} 0 & - \curl \\ \curl_0 & 0 \end{pmatrix}
\] is well-posed in $L^2_\eta(\R;L^2(\Omega)^6)$ for $\eta>0$ large enough.

\emph{Generalisations.} The material law operators introduced here are potentially more general than visible in the examples above. In fact, delay equations or equations with memory effects can be covered here, as well. Also the complexity of equations that can be dealt with under the explained approach are manifold. We refer to the literature for more examples, see e.g.~\cite{PTW15_WP_P}. For an extension to stochastic evolutionary equations we refer to \cite{SW16_SD} (autonomous case) and \cite{PTW18_SPDE} (non-autonomous case).

\section{Partial Differential Equations --\\ finite-dimensional nullspace}

This and the next section are concerned with convergence of evolutionary equations. Since we are aiming to cover homogenisation problems in particular, we are focussing on varying material laws. Thus, we are considering a sequence of evolutionary equations (or operators)
\[
   (M_n(\partial_t) +A)_n
\]
and want to address the limit $n\to\infty$. Similarly to the static case, we want to understand
\[
   (M_n(\partial_t) +A)U_n=F
\]
for fixed $F$. Thus, one is rather interested in the convergence of $(\overline{M_n(\partial_t)+A})^{-1}$. Given the results in the static case, we are again expecting convergence in the weak operator topology (although something more can be achieved in particular situations). In this section, we focus on the case of $A$ having compact resolvent. In the next section, we concentrate on $A=-A^*$ having compact resolvent, if reduced to the orthogonal complement of its nullspace.

To begin with, we analyse the topology of the material coefficients first.

\emph{A material law topology.} The topology, we endow the material laws with, is a topology of analytic functions such that pointwise, we use the weak-operator topology. It is thus natural to use a compact open topology (to keep the analyticity property) combined with the weak operator topology. The idea has emerged in \cite{W11_P,W12_HO}. For this let $E\subseteq \C$ be an open subset and $H,K$ be Hilbert spaces. Then we endow $\mathcal{H}(E;\mathcal{B}(H,K))$ with the initial topology so that
\[
      \mathcal{H}(E;\mathcal{B}(H,K))\ni M \mapsto (z\mapsto \langle \phi, M(z)\psi\rangle)\in \mathcal{H}(E)
\]
are continuous for all $\phi\in K$ and $\psi\in H$, where $\mathcal{H}(E)$ is endowed with the compact open topology, that is, uniform convergence on compact subsets of $E$. The resulting topological space is denoted by $\mathcal{H}_\w(E;\mathcal{B}(H,K))$.

\emph{A compactness statement.} The compactness of the compact open topology and the same for bounded subsets fo bounded linear operators under the weak operator topology leads to a compactness statement for the material law topology. Let $H,K$ be Hilbert spaces.
\begin{theorem}[{{{\cite[Theorem 4.3]{W14_FE}}}}] Let $\mathcal{N}\subseteq \mathcal{H}_\w(E;\mathcal{B}(H,K))$ be bounded. Then $\mathcal{N}$ is relatively compact. If both $H$ and $K$ are separable, then $\mathcal{N}$ is metrisable and $\overline{\mathcal{N}}$ is sequentially compact. 
\end{theorem}
\begin{corollary}[{{\cite[Theorem 3.4]{W12_HO}}}] Let $(M_n)_n$ in $\mathcal{H}_\w(E;\mathcal{B}(H,K))$ be bounded, $H,K$ separable. Then there exists $M\in \mathcal{H}_\w(E;\mathcal{B}(H,K))$  and $M_{n_k}\to M$ as $k\to\infty$ in $\mathcal{H}_\w(E;\mathcal{B}(H,K))$.
\end{corollary}
\emph{Convergence of evolutionary equations.} The convergence result that is underlying all classical homogenisation theorems of dynamic equations of mathematical physics is the following.
\begin{theorem}[{{\cite[Theorem 4.1]{W14_FE}}}]\label{thm:convcomp} Let $\mu>\nu\in \R$. Let $(M_n)_n$ be a convergent sequence in $\mathcal{H}_\w(\C_{\Re>\nu};\mathcal{B}(H,K))$; denote by $M$ its limit, $H,K$ Hilbert spaces. Assume that $\Re M_n(z)\geq \alpha$ for all $z\in \C_{\Re>\nu}$ and $n\in \N$ and assume that $z\mapsto z^{-1}M_n(z)$ is bounded uniformly in $n$. Assume that $A=-A^*$ with $\dom(A) \hookrightarrow\hookrightarrow H$. Then 
\[
     \overline{(M_n(\partial_t)+A)}^{-1} \to   \overline{(M(\partial_t)+A)}^{-1} \in \mathcal{B}_\w (L^2_\mu(\R;H)).
\]
\end{theorem}
Note that the compactness requirement is essential for the theorem to be true. In fact, the case is entirely different, if  $A=0$ and $H$ is infinite-dimensional. For this, we refer to the extensive studies \cite{W14_G} and \cite[Chapter 4]{W16_H} on \emph{ordinary} differential equations.

\emph{One-(plus-one)-dimensional example.} A prototype situation, where the compactness statement is easily verified is the case of one-dimensional heat conduction on $\Omega=(a,b)$ for some $a,b\in \R$. As above, we shall assume homogeneous Dirichlet boundary conditions. Thus, $A=-A^*=\begin{pmatrix} 0 & \partial_1 \\ \partial_{1,0} & 0 \end{pmatrix}$ has compact resolvent. Next, we let $a\in L^\infty(\R)$ be $1$-periodic with $\Re a\geq \alpha$; put $a_n \coloneqq a(n\cdot)$. Then by Theorem \ref{thm:convcomp}, we obtain
\begin{multline*}
    \overline{\partial_t \begin{pmatrix} 1 & 0 \\ 0 & 0 \end{pmatrix}+ \begin{pmatrix} 0 & 0 \\ 0 & a_n^{-1} \end{pmatrix} +\begin{pmatrix} 0 & \partial_1 \\ \partial_{1,0} & 0 \end{pmatrix}}^{-1}\\
    \to \overline{\partial_t \begin{pmatrix} 1 & 0 \\ 0 & 0 \end{pmatrix}+ \begin{pmatrix} 0 & 0 \\ 0 & \mathfrak{M}(a^{-1}) \end{pmatrix} +\begin{pmatrix} 0 & \partial_1 \\ \partial_{1,0} & 0 \end{pmatrix}}^{-1}\in  \mathcal{B}_\w (L^2_\nu(\R;(L^2(a,b))^2))
\end{multline*}for all $\nu>0$.
Note that the second order formulation of the limit equation would read
\[
    \partial_t \theta - \frac1{ \mathfrak{M}(a^{-1})}\partial_1\partial_{1,0}\theta =Q
\]for some right-hand side $Q$.

\emph{3-dimensional wave equation.} Let $\Omega\subseteq \R^3$ be open and bounded. In this example, we consider the convergence of the wave equation. For this let $(a_n)_n$ in $M_{\sym}(\alpha,\beta,\Omega)$ be $G$-converging to some $a$. In contrast to the previous example, we cannot directly apply the convergence statement to the operator
\[
     \partial_t \begin{pmatrix} 1 & 0 \\ 0 & a_n^{-1} \end{pmatrix} +\begin{pmatrix} 0 & \dive \\ \grad_0 & 0 \end{pmatrix}.
\]The reason for this is that $\dive$ has an \emph{inifinite-dimensional} nullspace. In order to rectify the situation, we revisit the second order formulation of the wave equation, which for some right-hand side $F$ reads
\[
    \partial_t^2 u -\dive a_n \grad_0 u = F.
\]We emphasise that we have used $0$ as an index to be reminded of Dirichlet-boundary conditions. Note that the second order formulation does not change, if we introduce the projection $\pi_{g_0}$ onto $g_0(\Omega)$ to the left of $\grad_0$ and to the right of $\dive$. Indeed, this is because we trivially have $\grad_0 u \in g_0(\Omega)$ and $\kar(\dive)=g_0(\Omega)^\bot$; hence $\dive = \dive \pi_{g_0}$. Using the canonical embedding $\iota_{g_0} \colon g_0(\Omega)\hookrightarrow L^2(\Omega)^3$, we get $\pi_{g_0} = \iota_{g_0}\iota_{g_0}^*$ so that the second order formulation reads
\[
  \partial_t^2 u -\dive \iota_{g_0}\iota_{g_0}^* a_n \iota_{g_0} \iota_{g_0}^* \grad_0 u = F.
\]Using $v=\partial_tu$ and $\tilde{\sigma}\coloneqq \iota_{g_0}^* a_n \iota_{g_0} \iota_{g_0}^* \grad_0 u$, we obtain
\[
\left(    \partial_t \begin{pmatrix} 1 & 0 \\ 0 & \left(\iota_{g_0}^* a_n \iota_{g_0}\right)^{-1} \end{pmatrix} +\begin{pmatrix} 0 & \dive\iota_{g_0} \\ \iota_{g_0}^*\grad_0 & 0 \end{pmatrix}\right)\begin{pmatrix} v \\ \tilde{\sigma} \end{pmatrix} = \begin{pmatrix} F \\ 0\end{pmatrix}.
\]
It is elementary to see that the latter equation is well-posed, by Theorem \ref{thm:st}. More importantly, we obtain that $\tilde{A}=\begin{pmatrix} 0 & \dive\iota_{g_0} \\ \iota_{g_0}^*\grad_0 & 0 \end{pmatrix}$ is skew-self-adjoint (see \cite[Lemma 4.4]{EGW17_D2N}) and has compact resolvent, see e.g.~\cite[Lemma 4.1 (also cf.~middle of p.~288)]{W13_HP}. Using the characterisation of $G$-convergence from Theorem \ref{thm:GH}(a), we finally obtain with Theorem \ref{thm:convcomp}
\begin{multline*}
  \overline{ \partial_t \begin{pmatrix} 1 & 0 \\ 0 & \left(\iota_{g_0}^* a_n \iota_{g_0}\right)^{-1} \end{pmatrix} +\begin{pmatrix} 0 & \dive\iota_{g_0} \\ \iota_{g_0}^*\grad_0 & 0 \end{pmatrix}}^{-1} \\
  \to   \overline{ \partial_t \begin{pmatrix} 1 & 0 \\ 0 & \left(\iota_{g_0}^* a \iota_{g_0}\right)^{-1} \end{pmatrix} +\begin{pmatrix} 0 & \dive\iota_{g_0} \\ \iota_{g_0}^*\grad_0 & 0 \end{pmatrix}}^{-1} \in \mathcal{B}_\w(L^2_\mu(\R;L^2(\Omega)\oplus g_0(\Omega)))
\end{multline*}for all $\mu>0$. Note that as a consequence of the latter result, we deduce that the solutions $u_n$ of 
\[
 \partial_t^2 u_n -\dive a_n \grad_0 u_n = F
\]
weakly converge (tested with bounded functions with compact support in time-space) to $u$, which solves
\[
   \partial_t^2 u -\dive a \grad_0 u = F.
\]

We emphasise that the `trick' of projecting away from the kernel of the divergence cannot be used in this form for Maxwell's equations. The reason for this is that both spatial derivative operators in Maxwell's equations consist of infinite-dimensional nullspaces. Moreover, there is no natural second order formulation, that would allow reducing the coefficients in question to certain range spaces. That is why, we need another step of generalisation, which will also help us to improve the result for the wave or heat equation.

\emph{Generalisations.} Although being of limited use for the full Maxwell system, the convergence Theorem \ref{thm:convcomp} has applications in one-dimensional transport problems, sub- and super fractional diffusion, coupled systems of the wave equation with an ordinary differential equation, and  thermo-elasticity: for all these examples see \cite[Section 4]{W13_HP}. Furthermore, applications to fractional elasticity have been found in \cite{W14_FE}; applications to highly oscillatory mixed type problems have been discussed in \cite{W16_SH,FW17_1D,Franz2018}. Quantitative results (optimal in order, operator-norm estimates) in this line of problems using compactness of the spatial derivative operator or a spectral gap type condition can be found in \cite{CW17_FH} (in arbitrary spatial dimensions) and in \cite{CW17_1D} (in a one-dimensional dynamic setting). Note that for quantitative estimates to be shown, in the mentioned results the periodicity of the coefficients is essential.

Furthermore, the whole theory admits an extension to non-autonomous partial differential equations; we refer the reader to the contribution \cite{Waurick2017} with applications to a homogenisation problem for a non-autonomous Kelvin-Voigt model, to homogenisation for acoustic waves problems with impedance type boundary conditions and a singular perturbation problem for a mixed type equation. A round up presentation of the results for non-autonomous equations can be found in \cite[Chapter 5]{W16_H}.

\section{Partial Differential Equations -- infinite-dimensional nullspace}

\emph{A decomposition of evolutionary equations.} The aim of this section is to understand the limit behaviour as $n\to\infty$ of 
\[
   (M_n(\partial_t) +A)_n,
\]
where in contrast to the previous section, we shall replace the compact embedding assumption of $\dom(A)$ into $H$ by asking for $\dom(A)\cap \kar(A)^\bot$ being compactly embedded into $H$, only. In order to apply the previous results, as well, we decompose $M_n(\partial_t)+A$ into a two-by-two block operator matrix acting on the direct sum $\rge(A)\oplus \kar(A)$. We will transform the resulting equation, as well. Thus, we need to take the right-hand side into account. We drop the index $n$ for a moment.
Consider
\[
   (M(\partial_t) +A)U=F.
\]
We obtain
\[
\left(     \begin{pmatrix} M(\partial_t)_{rr} & M(\partial_t)_{rk} \\ 
     M(\partial_t)_{kr} & M(\partial_t)_{kk}  \end{pmatrix} +\begin{pmatrix} A_{rr} & 0 \\ 0 & 0 \end{pmatrix}\right) \begin{pmatrix} U_r \\ U_k \end{pmatrix} = \begin{pmatrix} F_r \\ F_k \end{pmatrix},
\]
the index $r$ stands for $\rge(A)$ and $k$ for $\kar(A)$. The opertator $T_{rk}$ is the part mapped by $T$ from $\kar(A)$ into $\rge(A)$ (analogously, for  $T_{rr},T_{kr},T_{kk}$).
Multiplying this equation by $\begin{pmatrix}1 & - M(\partial_t)_{rk}M(\partial_t)_{kk} ^{-1} \\ 0 & M(\partial_t)_{kk}^{-1} \end{pmatrix}$, we obtain
\begin{multline*}
\left(     \begin{pmatrix} M(\partial_t)_{rr}- M(\partial_t)_{rk}M(\partial_t)_{kk} ^{-1}M(\partial_t)_{kr} & 0 \\ 
    M(\partial_t)_{kk}^{-1} M(\partial_t)_{kr} & 1  \end{pmatrix} +\begin{pmatrix} A_{rr} & 0 \\ 0 & 0 \end{pmatrix}\right) \begin{pmatrix} U_r \\ U_k \end{pmatrix} \\ = \begin{pmatrix} F_r - M(\partial_t)_{rk}M(\partial_t)_{kk} ^{-1}F_k \\ M(\partial_t)_{kk} ^{-1} F_k \end{pmatrix}.
\end{multline*}
A slight reformulation of this now is
\[
\begin{pmatrix}U_r \\ U_k \end{pmatrix} =    \begin{pmatrix} (M(\partial_t)_{rr}- M(\partial_t)_{rk}M(\partial_t)_{kk} ^{-1}M(\partial_t)_{kr}+A_{rr})^{-1}(F_r - M(\partial_t)_{rk}M(\partial_t)_{kk} ^{-1}F_k )  \\ 
   - M(\partial_t)_{kk}^{-1} M(\partial_t)_{kr} U_r + M(\partial_t)_{kk}^{-1} F_k \end{pmatrix}.
\]
\emph{Convergence statement.} The adapted convergence statement for the present situation requires an amended statement in Theorem \ref{thm:convcomp}, which reads as follows.
\begin{theorem}[{{\cite[Theorem 1.1]{W18_ONC_PAMM}}}]\label{thm:convcomp2} Let $\mu>\nu\in \R$, $H$ Hilbert space. Let $(M_n)_n$ a sequence in $\mathcal{H}_\w(\C_{\Re>\nu};\mathcal{B}(H))$ convergent to some $M.$  Assume that $\Re M_n(z)\geq \alpha$ for all $z\in \C_{\Re>\mu}$ and $n\in \N$ and assume that $z\mapsto z^{-1}M_n(z)$ is bounded uniformly in $n$. Assume that $A=-A^*$ with $\dom(A) \hookrightarrow\hookrightarrow H$. Then 
\[
     \partial_t^{-1}\overline{(M_n(\partial_t)+A)}^{-1} \to  \partial_t^{-1} \overline{(M(\partial_t)+A)}^{-1} \in \mathcal{B} (L^2_\mu(\R;H)).
\]and for all weakly convergent $(F_n)_n$ in $L^2_\mu(\R;H)$, we have
\[
    \overline{(M_n(\partial_t)+A)}^{-1}F_n \rightharpoonup   \overline{(M(\partial_t)+A)}^{-1}F \in L^2_\mu(\R;H).
\]
\end{theorem}
We emphasise the convergence of the operators $  \partial_t^{-1}\overline{(M_n(\partial_t)+A)}^{-1}$ in \emph{operator norm}. With this result and the reformulation outlined above, we arrive at convergence statements for evolutionary equations of the more general type discussed here.
\begin{theorem}[{{\cite[Theorem 5.5]{W16_HPDE}}}]\label{thm:conv} Let $\mu>\nu\in \R$, $H$ Hilbert space. Let $(M_n)_n$ be a sequence in $\mathcal{H}_\w(\C_{\Re>\nu};\mathcal{B}(H))$. Assume that $\Re M_n(z)\geq \alpha$ for all $z\in \C_{\Re>\nu}$ and $n\in \N$ and assume that $z\mapsto z^{-1}M_n(z)$ is bounded uniformly in $n$. Assume that $A=-A^*$ with $\dom(A)\cap\kar(A)^\bot \hookrightarrow\hookrightarrow H$.

Assume there exists $M\in \mathcal{H}_\w(\C_{\Re>\nu};\mathcal{B}(H))$ such that
\begin{align*}
  M_n(\partial_t)_{kk}^{-1} & \to  M(\partial_t)_{kk}^{-1} \in \mathcal{H}_\w(\C_{\Re>\nu};\mathcal{B}(\kar(A)) \\
  M_n(\partial_t)_{rk}M_n(\partial_t)_{kk} ^{-1} &\to M(\partial_t)_{rk}M(\partial_t)_{kk} ^{-1} \in \mathcal{H}_\w(\C_{\Re>\nu};\mathcal{B}(\rge(A),\kar(A))\\
  M_n(\partial_t)_{kk}^{-1} M_n(\partial_t)_{kr} &\to M(\partial_t)_{kk}^{-1} M(\partial_t)_{kr} \in \mathcal{H}_\w(\C_{\Re>\nu};\mathcal{B}(\kar(A), \rge(A))\\
    M_n(\partial_t)_{rr}- M_n(\partial_t)_{rk}M_n(\partial_t)_{kk} ^{-1}M_n(\partial_t)_{kr} & \to  M(\partial_t)_{rr}- M(\partial_t)_{rk}M(\partial_t)_{kk} ^{-1}M(\partial_t)_{kr} \\ &\quad\quad\quad\quad \in \mathcal{H}_\w(\C_{\Re>\nu};\mathcal{B}(\rge(A)).  
\end{align*}
Then
\[
   \overline{(M_n(\partial_t)+A)}^{-1} \to  \overline{(M(\partial_t)+A)}^{-1} \in \mathcal{B}_\w(L^2_\mu(\R;H)).
\]
\end{theorem}
\begin{remark} (a) In the case that $\kar(A)$ is infinite-dimensional, the convergence statement cannot be expected to be improved. This is already visible for $A=0$, see e.g.~\cite{W14_G}.

(b) Using the techniques of \cite[Proof of Theorem 5.5]{W18_NHC}, one can show that such an $M$ always exists (at least for a subsequence $(n_k)_k$).

(c) Note that \cite[Theorem 5.5]{W16_HPDE} as it is formulated in \cite{W16_HPDE} uses stronger assumptions. These, however, can be dispensed with entirely. The sketch of the proof is outlined above and rests on the decomposition of the evolutionary equation in range and kernel space of $A$ and some suitable positive definiteness estimates.
\end{remark}
The convergence result in Theorem \ref{thm:conv} is at the heart of nonlocal $H$-convergence; particularly, when it comes to Maxwell's equations. Before, however, we finalise this paper with applications to electro-magentics, we revisit the $3$-dimensional heat equation.

\emph{First order heat conduction.}  Let $\Omega\subseteq \R^3$ be a standard domain. Assume that $(a_n)_n$ in $M(\alpha,\beta,\Omega)$ is $H$-convergent to some $a$. We address the convergence of the operator sequence
\[
    \left( \overline{ \partial_t \begin{pmatrix} 1 & 0 \\ 0 & 0 \end{pmatrix}+ \begin{pmatrix} 0 & 0 \\ 0 & a_n^{-1} \end{pmatrix} +\begin{pmatrix} 0 & \dive \\ \grad_0 & 0 \end{pmatrix}}^{-1}\right)_n
\]in $\mathcal{B}_\w(L^2(\R;L^2(\Omega)^{3+1}))$. For this, we need to analyse, whether the conditions in Theorem \ref{thm:conv} are met. Let us have a closer look at the range and the kernel of $A=\begin{pmatrix} 0 & \dive \\ \grad_0 & 0 \end{pmatrix}$. We have
\begin{align*}
    \rge(A)&=\rge(\dive)\oplus \rge(\grad_0)=L^2(\Omega)\oplus g_0(\Omega)\, \text{and} \\ \kar(A)&=\kar(\grad_0)\oplus \kar(\dive)=\{0\}\oplus c(\Omega),
\end{align*}
where we have used Theorem \ref{thm:HD}. Due to the block structure of $M_n(\partial_t)= \partial_t \begin{pmatrix} 1 & 0 \\ 0 & 0 \end{pmatrix}+ \begin{pmatrix} 0 & 0 \\ 0 & a_n^{-1} \end{pmatrix} $, we only need to decompose $a_n^{-1}$ according to the (second component of the) range and kernel space of $A$, that is, along the decomposition $g_0(\Omega)\oplus c(\Omega)$. Using $\iota_{g_0}$ and $\iota_c$ to be the canonical embeddings, we thus need to confirm the convergence of
\begin{align*}
  &   \left(\iota_c^* a_n^{-1} \iota_c\right)^{-1} \to    \left( \iota_c^* a^{-1} \iota_c\right)^{-1} \\
  &  \iota_{g_0} a_n^{-1} \iota_c \left(\iota_c^* a_n^{-1} \iota_c\right)^{-1} \to   \iota_{g_0} a^{-1}\iota_c \left( \iota_c^* a^{-1} \iota_c\right)^{-1} \\
    &  \left(\iota_c^* a_n^{-1} \iota_c\right)^{-1}\iota_{c}^* a_n^{-1} \iota_{g_0}  \to  \left( \iota_c^* a^{-1} \iota_c\right)^{-1}\iota_{c}^* a^{-1} \iota_{g_0} \\
    & \iota_{g_0}^* a_n^{-1} \iota_{g_0} - \iota_{g_0}^* a_n^{-1} \iota_{c}  \left(\iota_c^* a_n^{-1} \iota_c\right)^{-1}\iota_{c}^* a_n^{-1} \iota_{g_0}  \to \iota_{g_0}^* a^{-1} \iota_{g_0} - \iota_{g_0}^* a^{-1} \iota_{c}  \left(\iota_c^* a^{-1} \iota_c\right)^{-1}\iota_{c}^* a^{-1} \iota_{g_0}
\end{align*}in the respective weak operator topologies. With the representation as block operator matrices, the latter expressions can be simplified. In fact, we have
\begin{align*}
  &   \left(\iota_c^* a_n^{-1} \iota_c\right)^{-1} = a_{n,11}-a_{n,10}a_{n,00}^{-1}a_{n,01} \\
  &  \iota_{g_0} a_n^{-1} \iota_c \left(\iota_c^* a_n^{-1} \iota_c\right)^{-1} = a_{n,00}^{-1}a_{n,01} \\
    &  \left(\iota_c^* a_n^{-1} \iota_c\right)^{-1}\iota_{c}^* a_n^{-1} \iota_{g_0} = a_{n,10}a_{n,00}^{-1} \\
    & \iota_{g_0}^* a_n^{-1} \iota_{g_0} - \iota_{g_0}^* a_n^{-1} \iota_{c}  \left(\iota_c^* a_n^{-1} \iota_c\right)^{-1}\iota_{c}^* a_n^{-1} \iota_{g_0} = a_{n,00}^{-1},
   \end{align*}where we employed the notation from the definition of nonlocal $H$-convergence. Note that the equalities can be seen most easily with the Schur complement formulas outlined in \cite[Lemma 4.8]{W18_NHC}.
   Thus, as (local) $H$-convergence implies nonlocal $H$-convergence (Theorem \ref{thm:GH}(b)), we finally arrive at the following convergence result
   \begin{multline*}
    \overline{ \partial_t \begin{pmatrix} 1 & 0 \\ 0 & 0 \end{pmatrix}+ \begin{pmatrix} 0 & 0 \\ 0 & a_n^{-1} \end{pmatrix} +\begin{pmatrix} 0 & \dive \\ \grad_0 & 0 \end{pmatrix}}^{-1} \\
         \to    \overline{ \partial_t \begin{pmatrix} 1 & 0 \\ 0 & 0 \end{pmatrix}+ \begin{pmatrix} 0 & 0 \\ 0 & a^{-1} \end{pmatrix} +\begin{pmatrix} 0 & \dive \\ \grad_0 & 0 \end{pmatrix}}^{-1} \in \mathcal{B}_\w(L^2_\nu(\R;L^2(\Omega)^{3+1}))
   \end{multline*}as $n\to\infty$. The upshot of this result is that not only the heat but also the heat flux weakly converges in the limit. This is an implication that has not been available with Theorem \ref{thm:convcomp}, which has been exemplified with the wave equation above.

\emph{Maxwell's equations.} Finally, we shall treat an example with nonlocal $H$-convergence, see \cite[Section 7]{W18_NHC} for more details. As we want to illustrate the results rather than providing a sophisticated example, we assume that the electric conductivity vanishes, that is, $\sigma=0$. Let $\Omega\subseteq \R^3$ be a standard domain. Let $(\varepsilon_n)_n$ be a convergent sequence in $({\mathcal{M}}(\alpha,\beta,\Omega),{\tau}_{\nlH})$ and $(\mu_n)_n$ convergent sequence in $(\hat{\mathcal{M}}(\alpha,\beta,\Omega),\hat{\tau}_{\nlH})$. Assuming that $\varepsilon_n=\varepsilon_n^*, \mu_n=\mu_n^*\geq \alpha$ for all $n\in\N$, we obtain
\[
   \overline{   \partial_t \begin{pmatrix} \varepsilon_n & 0 \\ 0 & \mu_n \end{pmatrix} + \begin{pmatrix} 0 & - \curl \\ \curl_0 & 0 \end{pmatrix}}^{-1}
    \to     \overline{   \partial_t \begin{pmatrix} \varepsilon & 0 \\ 0 & \mu \end{pmatrix} + \begin{pmatrix} 0 & - \curl \\ \curl_0 & 0 \end{pmatrix}}^{-1}\in \mathcal{B}_\w(L^2_\nu(\R;L^2(\Omega)^6))
\] for all $\nu>0$. 

We shall argue next, why the latter convergence statement indeed follows from Theorem \ref{thm:conv}. To begin with, we focus on the compactness condition imposed on $A$. This is proven by using Picard's selection theorem:
\begin{theorem}[\cite{Picard1984}] Let $\Omega\subseteq \R^3$ be a bounded, weak Lipschitz domain. Then $\dom(\curl)\cap\dom(\dive_0)\hookrightarrow\hookrightarrow L^2(\Omega)^3$ and  $\dom(\curl_0)\cap\dom(\dive)\hookrightarrow\hookrightarrow L^2(\Omega)^3$.
\end{theorem}

The convergence assumptions in Theorem \ref{thm:conv} are precisely the ones yielded by the assumption on nonlocal $H$-convergence; see also \cite[Section 7]{W18_NHC} (note the misprint in Example 7.10; the roles of $\varepsilon$ and $\mu$ need to be swapped in the second to last and in the last line). Indeed, since $A=\begin{pmatrix} 0 & -\curl \\ \curl_0 & 0 \end{pmatrix}$, we obtain that 
\begin{align*}
     \rge(A) & = c(\Omega) \oplus c_0(\Omega) \\
     \kar(A) & = g_0(\Omega) \oplus g(\Omega).
\end{align*}
Hence, with $M_n(\partial_t) =   \partial_t \begin{pmatrix} \varepsilon_n & 0 \\ 0 & \mu_n \end{pmatrix}$, we obtain for all $n\in\N$
\begin{align*}
& z  M_n(z)_{kk}^{-1}  = \begin{pmatrix} (\iota_{g_0}^* \epsilon_n \iota_{g_0})^{-1} & 0 \\ 0 & (\iota_{g}^* \mu_n \iota_{g})^{-1} \end{pmatrix} \\
  & M_n(z)_{rk}M_n(z)_{kk} ^{-1}  =  \begin{pmatrix} (\iota_{c}^* \epsilon_n \iota_{g_0}) & 0 \\ 0 &(\iota_{c_0}^* \mu_n \iota_{g}) \end{pmatrix} \begin{pmatrix} (\iota_{g_0}^* \epsilon_n \iota_{g_0})^{-1} & 0 \\ 0 & (\iota_{g}^* \mu_n \iota_{g})^{-1} \end{pmatrix} \\
  &    M_n(z)_{kk} ^{-1}M_n(z)_{kr}  =  \begin{pmatrix} (\iota_{g_0}^* \epsilon_n \iota_{g_0})^{-1} & 0 \\ 0 & (\iota_{g}^* \mu_n \iota_{g})^{-1} \end{pmatrix} \begin{pmatrix} (\iota_{g_0}^* \epsilon_n \iota_{c}) & 0 \\ 0 & (\iota_{g}^* \mu_n \iota_{c_0}) \end{pmatrix} \\
 &  z^{-1}  M_n(z)_{rr}- z^{-1}M_n(z)_{rk}M_n(z)_{kk} ^{-1}M_n(z)_{kr}  =  \begin{pmatrix} \iota_{c}^* \epsilon_n \iota_{c} & 0 \\ 0 &  \iota_{c_0}^* \mu_n \iota_{c_0} \end{pmatrix}\\ & \quad\quad\quad\quad\quad\quad\quad-  \begin{pmatrix} (\iota_{c}^* \epsilon_n \iota_{g_0}) & 0 \\ 0& (\iota_{c_0}^* \mu_n \iota_{g}) \end{pmatrix} \begin{pmatrix} (\iota_{g_0}^* \epsilon_n \iota_{g_0})^{-1} & 0 \\ 0 & (\iota_{g}^* \mu_n \iota_{g})^{-1} \end{pmatrix} \begin{pmatrix} (\iota_{g_0}^* \epsilon_n \iota_{c}) & 0 \\ 0 & (\iota_{g}^* \mu_n \iota_{c_0}) \end{pmatrix}.
\end{align*}
Using the notation from the definition of nonlocal $H$-convergence, we can reformulate the right-hand sides to obtain
\begin{align*}
& \begin{pmatrix} (\iota_{g_0}^* \epsilon_n \iota_{g_0})^{-1} & 0 \\ 0 & (\iota_{g}^* \mu_n \iota_{g})^{-1} \end{pmatrix} = \begin{pmatrix} \epsilon_{n,00}^{-1} & 0 \\ 0 & \hat{\mu}_{n,00}^{-1} \end{pmatrix}\\
&  \begin{pmatrix} (\iota_{c}^* \epsilon_n \iota_{g_0}) & 0 \\ 0 &(\iota_{c_0}^* \mu_n \iota_{g}) \end{pmatrix} \begin{pmatrix} (\iota_{g_0}^* \epsilon_n \iota_{g_0})^{-1} & 0 \\ 0 & (\iota_{g}^* \mu_n \iota_{g})^{-1} \end{pmatrix} =  \begin{pmatrix}  \epsilon_{n,10}\epsilon_{n,00}^{-1} & 0 \\ 0 &\hat{\mu}_{n,10} \hat{\mu}_{n,00}^{-1} \end{pmatrix}\\
  &  \begin{pmatrix} (\iota_{g_0}^* \epsilon_n \iota_{g_0})^{-1} & 0 \\ 0 & (\iota_{g}^* \mu_n \iota_{g})^{-1} \end{pmatrix} \begin{pmatrix} (\iota_{g_0}^* \epsilon_n \iota_{c}) & 0 \\ 0 & (\iota_{g}^* \mu_n \iota_{c_0}) \end{pmatrix} =
  \begin{pmatrix}  \epsilon_{n,00}^{-1}\epsilon_{n,01} & 0 \\ 0 & \hat{\mu}_{n,00}^{-1}\hat{\mu}_{n,01} \end{pmatrix} \\
 &   \begin{pmatrix} \iota_{c}^* \epsilon_n \iota_{c} & 0 \\ 0 &  \iota_{c_0}^* \mu_n \iota_{c_0} \end{pmatrix}-  \begin{pmatrix} (\iota_{c}^* \epsilon_n \iota_{g_0}) & 0 \\ 0& (\iota_{c_0}^* \mu_n \iota_{g}) \end{pmatrix} \begin{pmatrix} (\iota_{g_0}^* \epsilon_n \iota_{g_0})^{-1} & 0 \\ 0 & (\iota_{g}^* \mu_n \iota_{g})^{-1} \end{pmatrix} \begin{pmatrix} (\iota_{g_0}^* \epsilon_n \iota_{c}) & 0 \\ 0 & (\iota_{g}^* \mu_n \iota_{c_0}) \end{pmatrix} \\ &=   \begin{pmatrix} \epsilon_{n,11}- \epsilon_{n,10}\epsilon_{n,00}^{-1}\epsilon_{n,01} & 0 \\ 0 & \hat{\mu}_{n,11}-\hat{\mu}_{n,10}\hat{\mu}_{n,00}^{-1}\hat{\mu}_{n,01} \end{pmatrix},
\end{align*}
which are precisely the quantities we have assumed to be convergent in the weak operator topology as $(\epsilon_n)_n$ and $(\mu_n)_n$ converge to $\epsilon$ and $\mu$ in $\tau_{\nlH}$ and $\hat{\tau}_{\nlH}$, respectively.

\emph{Generalisations.} The main theorem of this section (Theorem \ref{thm:conv}) has been applied to bi-anisotropic dissipative media in electro-magnetics and to thermo-piezo-electricity, see \cite[Section 6]{W16_HPDE}. The Theorem \ref{thm:conv}, however, is written in a way that it applies to even more complicated autonomous equations and systems of mathematical physics. 

\section{Concluding remarks}

In this surveying article, we have highlighted the intimate relationship of homogenisation problems and the weak operator topology. We have highlighted characterisations in the static case and have shown several implications for time-dynamic problems. The strategies developed form a profound basis for all classical homogenisation theorems for linear problems of mathematical physics. Non-classical problems with nonlocal coefficients can be treated as well.

We conclude with some unresolved issues. Bounded subsets of the set bounded linear operators are metrisable under the weak operator topology. A main line of study would thus be a convergence rate analysis of the dynamic problem given precise rates for the static variants are known. Convergence rates for evolutionary equations of a particular type are known, though. These cases, however, are restricted to periodicity assumptions on the coefficients. If the coefficients are not uniformly (in $n$) bounded away from zero, the techniques outlined above are not applicable anymore. This, however, is the main issue of so-called high-contrast homogenisation problems, \cite{Cherednichenko2016}. There have been several contributions to non-autonomous equations, as well, see \cite{W16_H,Waurick2017}. However, the case of non-autonomous Maxwell's equations is still open. Coming back to the static case, it is ongoing research to understand the concept of nonlocal $H$-convergence also in the case of non-standard domains $\Omega$.

\bibliographystyle{abbrv}

\noindent
Marcus Waurick \\Department of Mathematics and Statistics, University of Strathclyde,\\
Livingstone Tower, 26 Richmond Street, Glasgow G1 1XH,\\
Scotland\\
Email: %M.Waurick@bath.ac.uk
{\tt marcus.wau\rlap{\textcolor{white}{hugo@egon}}rick@strath.\rlap{\textcolor{white}{darmstadt}}ac.uk}

\end{document}